\documentclass[a4paper,12pt]{article}
\usepackage{amssymb}
\usepackage{amsfonts}
\usepackage{amsmath}
\usepackage{bbm}

\newcommand{\slnc}{\mathfrak{sl}(N,\mathbb{C})}
\newcommand{\C}{\mathbb{C}}

\newcommand{\np}{\mathfrak{n}_+}
\newcommand{\nm}{\mathfrak{n}_-}
\newcommand{\g}{\mathfrak{g}}
\newcommand{\h}{\mathfrak{h}}
\newcommand{\ap}{\mathfrak{a}}
\newcommand{\am}{\widetilde{\mathfrak{a}}}
\newcommand{\bp}{\mathfrak{b}}
\newcommand{\bm}{\widetilde{\mathfrak{b}}}
\newcommand{\up}{U}
\newcommand{\um}{\widetilde{U}}
\newcommand{\vp}{V}
\newcommand{\vm}{\widetilde{V}}
\newcommand{\indp}{J_k}
\newcommand{\indm}{\widetilde{J}_k}
\newcommand{\ad}[1]{\mathrm{ad}_{#1}}
\newcommand{\proof}{{\noindent\bf Proof: }}
\newcommand{\qed}{\hfill$\Box$}
\newcommand{\ih}{i_\mathfrak{h}}

\DeclareMathOperator{\spn}{span}

\newtheorem{Theorem}{Theorem}[section]

\newtheorem{Lemma}[Theorem]{Lemma}
\newtheorem{Corollary}[Theorem]{Corollary}

\begin{document}
\title{Wei-Norman equations for a unitary evolution}

\author{Szymon Charzy\'nski$^1$ and Marek Ku\'{s}$^2$
\\
 \\
$^1$Faculty of Mathematics and Natural Sciences, \\ Cardinal Stefan
Wyszy\'nski University, \\ ul. W\'oycickiego 1/3, 01-938 Warszawa, Poland
\\
\\
 $^2$Center for Theoretical Physics, Polish Academy of Sciences, \\
 Aleja Lotnik\'ow 32/46, 02-668 Warszawa, Poland \\
} \maketitle
\begin{abstract}
The Wei-Norman technique allows to express the solution of a system of linear
non-autonomous differential equations in terms of product of exponentials. In
particular it enables to find a time-ordered product of exponentials by solving
a set of nonlinear differential equations. The method has numerous theoretical
and computational advantages, in particular in optimal control theory. We show
that in the unitary case, i.e.\ when the solution of the linear system is given
by a unitary evolution operator, the nonlinear system can be by an appropriate
choice of ordering reduced to a hierarchy of matrix Riccati equations. Our
findings have a particular significance in quantum control theory since pure
quantum evolution is unitary.
\end{abstract}

\section{Introduction}

In two papers \cite{wei63} and \cite{wei64} written nearly fifty years ago,
Wei and Norman developed a method for solving systems of linear differential
equations with variable coefficients based on Lie-group techniques. As it is
well known, solutions of non-autonomous linear systems of differential
equations can be expressed in terms of time ordered exponentials. From this
point of view one can treat Wei-Norman formulae as a way to calculate such
exponentials. More generally they provide local coordinates on a smooth
manifold on which acts a (finite-dimensional) Lie group of transformations in
terms of exponential map from the Lie algebra.

Since then the Wei-Norman method has found numerous applications in control and
system theory, see, e.g.\ \cite{brockett72}, \cite{brockett73},
\cite{huillet87}, \cite{leonard93}, \cite{chiou94}, \cite{leonard95},
\cite{charalambous00}, \cite{carinena03}, as well as a basis for numerical
approximation methods of integrating linear non-autonomous systems
\cite{owren01}, \cite{celledoni01}, \cite{zanna02}.

In recently rapidly developing theory of quantum optimal control, (see e.g.\
\cite{carinena02}, \cite{carinena03a}, \cite{carinena09}, \cite{kuna10},
\cite{nihtila11}) of particular importance is the case when the underlying Lie
algebra is $\mathfrak{su}(N)$, since the pure quantum evolution preserves the
norm of a state vector. We show that in this case the Wei-Norman method leads
to a hierarchy of matrix Riccati equations, provided that a proper ordering of
generators (basis Lie algebra elements) is chosen. Our findings hinge
substantially on the structure of the $\mathfrak{su}(N)$ algebra and its
complexification $\mathfrak{sl}(N,\mathbb{C})$, in particular on the properties
of commutative ideals of the latter.

In the next section we present briefly the basics of the Wei-Norman method for
general complex semisimple Lie group. In Section~\ref{sec:slN} we summarize
properties of the $\mathfrak{sl}(N,\mathbb{C})$ relevant for our reasoning, in
particular the structure of its abelian ideals. Results of this section are
translated in Sections~\ref{sec:adjoint} and \ref{sec:nilpotent} to properties
of the adjoint endomorphism and exponential function on
$\mathfrak{sl}(N,\mathbb{C})$. The final result for the Wei-Norman formulae in
the case of $\mathfrak{sl}(N,\mathbb{C})$ is given in Section~\ref{sec:basis}.
The corresponding results for $\mathfrak{su}(N)$ is easy to recover by
restriction to anti-Hermitian generators. In Section~\ref{sec:examples} we give
explicit formulae for low dimensions.

\section{General Wei-Norman method}
\label{sec:WN}

Let $G$ be a $n$-dimensional Lie group and $\mathfrak{g}$ - its Lie algebra.
We assume in the following that $\mathfrak{g}$ is complex and simple. Let
also $\mathbb{R}\ni t\mapsto M(t)\in\mathfrak{g}$ be a curve in
$\mathfrak{g}$ and $K(t)$ - a curve in $G$ given by the differential equation
\begin{equation}\label{precession}
\frac{d}{dt}K(t)=M(t)K(t),\quad K(0)=I.
\end{equation}
In $\mathfrak{g}$ we choose some basis $X_k$, $k=1,\ldots,n$ in which $M(t)$
takes the form
\begin{equation}\label{Mdec}
M(t)=\sum_{k=1}^na_k(t)X_k,
\end{equation}
and look for the solution $K(t)$ in the form
\begin{equation}\label{coord2}
K(t)=\prod_{k=1}^n\exp\big(u_k(t)X_k\big),
\end{equation}
involving $n$ unknown functions $u_k$. Differentiating (\ref{coord2}) we get
\begin{eqnarray}\label{rach1}
K^\prime&=&
\sum_{l=1}^nu_l^\prime\prod_{k=1}^{l-1}\exp( u_kX_k)X_l
\prod_{k=l}^n\exp( u_kX_k) \nonumber \\
&=&\sum_{l=1}^n u_l^\prime\prod_{k=1}^{l-1}\exp( u_kX_k)X_l
\prod_{k=l-1}^1\exp(- u_kX_k)\prod_{k=1}^n\exp( u_kX_k) \nonumber \\
&=&\sum_{l=1}^n u_l^\prime\prod_{k=1}^{l-1}
\mathrm{Ad}_{\exp( u_kX_k)}\cdot X_l\,K
=\sum_{l=1}^n u_l^\prime\prod_{k<l}
\exp( u_k\mathrm{ad}_{X_k})\cdot X_l\,K,
\end{eqnarray}
where by $^\prime$ we denoted the differentiation with respect to $t$ and
$\mathrm{Ad}$ is the adjoint action of $G$ on $\mathfrak{g}$,
\begin{equation}\label{Ad}
\mathrm{Ad}_g\cdot X:= gXg^{-1}, \quad g\in G, \quad X\in\mathfrak{g}.
\end{equation}
In the last equality we used used $\mathrm{Ad}_{\exp( b X)}=\exp(
b\,\mathrm{ad}_X)$, where $\mathrm{ad}_X=[X,\cdot]$ is the adjoint action of
$\mathfrak{g}$ on itself.

Comparing (\ref{precession}) i (\ref{rach1}) we obtain
\begin{equation}\label{eq-wn0}
M(t)=\sum_{l=1}^nu_l^\prime\prod_{k<l}
\exp(u_k\mathrm{ad}_{X_k})\cdot X_l.
\end{equation}
Both sides of (\ref{eq-wn0}) are elements of $\mathfrak{g}$ and expanding
both in the basis $X_k$ we arrive at a system of coupled differential
equations for $u_k$ in terms of $a_l$, $k,l=1,\ldots,n$.

It is worthwhile to rewrite (\ref{eq-wn0}) in more compact form. Let us denote
\begin{equation}\label{}
A^{(l)}=\prod_{k<l}\exp( u_k\mathrm{ad}_{X_k}),
\end{equation}
and, consequently
\begin{equation}\label{Al}
A^{(l)}\cdot X_l=\sum_{j=1}^nA^{(l)}_{jl}X_j.
\end{equation}
Now (\ref{eq-wn0}) can be written as
\begin{equation}\label{eq-wn1}
\sum_{j=1}^n a_jX_j=\sum_{l=1}^n u_l^\prime\sum_{j=1}^n A^{(l)}_{jl} X_j
=\sum_{j=1}^n\left(\sum_{l=1}^n A^{(l)}_{jl} u_l^\prime\right)X_j,
\end{equation}
hence\
\begin{equation}\label{eq-wn2}
a_j=\sum_{l=1}^n A^{(l)}_{jl} u_l^\prime,
\end{equation}
or in a compact form
\begin{equation}\label{eq-wn2a}
\mathbf{a}=A\mathbf{u}^\prime,
\end{equation}
where $A$ is an $n\times n$ matrix with elements $A_{jl}=A^{l}_{jl}$, i.e., its
$l$-th column is equal to the $l$-th column of the matrix $A^{(l)}$, c.f.
(\ref{Al}), and $\mathbf{a}$ nad $\mathbf{u}$ are the vectors of the
coefficients of $M$ in (\ref{Mdec}) and the unknowns $u_k$. If $A$ is
invertible we obtain thus a system of (nonlinear) differential equations solved
for the first derivatives
\begin{equation}\label{eq-wn3}
\mathbf{u}^\prime=A^{-1}\mathbf{a}.
\end{equation}

\section{The structure of $\mathfrak{sl}(N,\mathbb{C})$ algebra}
\label{sec:slN}

Let us start with a brief summary of the structure of simple Lie algebras
\cite{humphreys}.

Each complex simple Lie algebra $\mathfrak{g}$ can be decomposed into the
root spaces with respect to a chosen Cartan subalgebra $\mathfrak{h}$ (a
maximal commutative subalgebra of $\mathfrak{g}$)
\begin{equation}\label{rootdecomp}
\mathfrak{g}=\mathfrak{h}\oplus\mathop{\bigoplus}\limits_{\alpha\in\Delta}
\mathfrak{g}_\alpha,
\end{equation}
where the one-dimensional \textit{root spaces} are defined as
\begin{equation}\label{rootspaces}
\mathfrak{g}_\alpha:=\{X\in\mathfrak{g}:[H,X]=\alpha(H)X\;
\forall H\in\mathfrak{h}\}
\end{equation}
The linear forms $\alpha\in\mathfrak{h}^*$ (the dual space to the algebra
$\mathfrak{h}$) are called \emph{roots} and the element $X_\alpha$ spanning the
subspace $\mathfrak{g}_\alpha$ is called the \emph{root vector}. By fixing a
basis $\{H_1,\ldots,H_r\}$ in $\mathfrak{h}$ we can identify the set of roots
$\Delta$ with a set of $r$-dimensional vectors
$\alpha=(\alpha^1,\ldots\alpha^r)$, where $\alpha^k=\alpha(H_k)$. The number
$r=\dim\mathfrak{h}$ is called the rank of $\mathfrak{g}$. Among the roots we
can find a set $\Pi$ of $r$ \textit{positive simple roots},
$\Pi=\{\alpha_1,\ldots\alpha_r\}$, such that for each $\alpha\in\Delta$
\begin{equation}\label{simpleroots}
\alpha=\sum_{i=1}^r m_i\alpha_i,
\end{equation}
where either all $m_i$ are nonnegative (such roots $\alpha$ for a set of
\textit{positive roots} $\Delta_+$) or all $m_i$ are nonpositive (such roots
$\alpha$ forming the set $\Delta_-$ of negative roots). There is one-to-one
correspondence between positive and negative roots: for each positive root
$\alpha$ the is a negative one $-\alpha$.

The root spaces have the following property
\begin{equation}\label{commg}
[\mathfrak{g}_\alpha,\mathfrak{g}_\beta]\subset\mathfrak{g}_{\alpha+\beta},
\end{equation}
(for $\beta=-\alpha$ we have
$[\mathfrak{g}_\alpha,\mathfrak{g}_{-\alpha}]\subset\mathfrak{h}$).

In terms of positive and negative roots the decomposition (\ref{rootdecomp})
can be rewritten as
\begin{equation}\label{cartan}
\mathfrak{g}=\mathfrak{n}_-\oplus\mathfrak{h}\oplus\mathfrak{n}_+, \quad
\mathfrak{n}_\pm=\mathop{\bigoplus}\limits_{\alpha\in\Delta_\pm}
\mathfrak{g}_\alpha.
\end{equation}
The subalgebras
\begin{equation}\label{borel}
\mathfrak{b_\pm}=\mathfrak{h}\oplus\mathfrak{n_\pm}
\end{equation}
are called the Borel subalgebras relative to the Cartan subalgebra
$\mathfrak{h}$.

In the following we take $G=SL(N, \mathbb{C})$ and $\g=\slnc$. In
this case $n=N^2-1$. To construct a basis in $\mathfrak{g}$ we
define
\begin{equation}\label{Skl}
S_{kl}=e_k e_l^\dagger,
\end{equation}
where $e_k$, $k=1,\ldots,N$, is the standard basis in $\mathbb{C}^N$, ie.,
$(e_k)_l=\delta_{kl}$. The commutation relations for $S_{kl}$ read
\begin{equation}\label{commut}
 [S_{ij},S_{kl}]=\delta_{kj}S_{il}-\delta_{il}S_{kj}.
\end{equation}
For $k\neq l$ elements $S_{pq}$ are the root vectors. The root
corresponding to roots space spanned by $S_{pq}$ will be denoted by
$\alpha_{pq}$. From $S_{kl}$ we construct a basis $X_m$ by the
following renumbering


\begin{eqnarray}
X_m&=&S_{pq}, \quad m=\frac{(N+q-1)(N-q)}{2}+p, \quad \mathrm{for\  } p<q, \label{Bu} \\
X_m&=&S_{ll}-S_{l+1,l+1}, \quad m=\frac{1}{2}N(N-1)+l, \quad l=1,\ldots,N-1, \label{Bcartan} \\
X_m&=&S_{qp}, \quad m=N^2-\frac{(N+q-1)(N-q)}{2}+p, \quad \mathrm{for\  } p<q.  \label{Bl}
\end{eqnarray}
One may choose roots $\alpha_{pq}$ corresponding to root vectors $S_{pq}$ for
$p<q$ to be the positive roots. Then the matrices (\ref{Bu}) and (\ref{Bl})
generate maximal nilpotent subalgebras $\mathfrak{n}_+$ and $\mathfrak{n}_-$
respectively  (conjugated with respect to the standard Hermitian structure on
$\mathbb{C}^N$), whereas (\ref{Bcartan}) generate the Cartan subalgebra of
$\slnc$. In this basis $\np$ is the set of strictly  upper triangular matrices
and $\nm$ is the set of strictly lower triangular matrices.

The set $\Pi$ of positive simple roots consist of elements of
$\alpha_{l,l+1}$ for $l=1,\ldots, N-1$.

In what follows we will show many useful features of this basis. First observe
that according to (\ref{Bu}), $\np:=\spn\{X_1,\ldots,X_{N(N-1)/2}\}$ and the
order of root vectors spanning $\np$ corresponds to the following order of
roots:
\begin{equation}\label{ord_np}
\alpha_{i,j}\succ\alpha_{k,l} \Leftrightarrow j<l \mbox{ or } (j=l
\mbox{ and } i>k).
\end{equation}

Analogously  $\nm:=\spn\{X_{N(N+1)/2},\ldots,X_{N^2-1}\}$ and the order of root
vectors spanning $\nm$ corresponds to the following order of roots:
\begin{equation}\label{ord_nm}
\alpha_{i,j}\succ\alpha_{k,l} \Leftrightarrow i>k \mbox{ or } (i=k
\mbox{ and } j<l).
\end{equation}

Recall, that by the theorem of Lie \cite{humphreys}, if $L$ is a nilpotent Lie
algebra of dimension $r$ then there exist a sequence of ideals $I_k$ fulfilling
\begin{equation}
0=I_0\subset I_1\subset \ldots \subset I_{r-1} \subset I_r =L
\end{equation}
and a basis $Y_1,\ldots,Y_r$ in $L$, such that
\begin{equation}
I_k:=\spn\{Y_1,\ldots,Y_k\}.
\end{equation}
The basis $\{X_1,\ldots,X_{N(N-1)/2}\}$ of $\np$ defined in (\ref{Bu}) has this
property. It follows from (\ref{commg}) and the fact that
$\alpha_{i,j}+\alpha_{k,l}\in\Delta$ provided $i=l$ or $j=k$. In both cases
$\alpha_{i,j}+\alpha_{k,l}\prec\alpha_{i,j}$, as can be easily checked. Thus
$I_k:=\spn\{X_1,\ldots,X_k\}$ is an ideal in $\np$ for any $k=1,\ldots,\frac 12
N(N-1)$.

Direct consequence of Lie theorem is that in this basis the matrix of
endomorphism $\ad{X_k}:\np\to\np$ is strictly upper triangular
\cite{humphreys}. Moreover, according to (\ref{Bl}) matrices generating $\nm$
are transpositions of matrices generating $\nm$ ordered in the reverse order.
This implies that matrices of endomorphisms $\ad{X_k}:\nm\to\nm$ for
$k=\frac12N(N+1),\ldots,N^2-1$ are strictly lower triangular.

We may decompose $\np$ in a particular direct sum of subspaces to
modify slightly the original Wei-Norman expansion (\ref{coord2}). To
this end let us consider the family $\{\mathfrak{a}_k\}$,
$k=1,\ldots,N-1$ of subspaces (the numbering of the root vectors
$X_i$ is the same as in (\ref{Bu})),
\begin{eqnarray}\label{setideals}
\mathfrak{a}_k=\spn\{X_{i_k},X_{i_k+1}\ldots,X_{i_k+N-k-1}\}, \quad
i_k=N(k-1)-\frac{k(k-1)}{2}+1.
  \end{eqnarray}
In the standard matrix representation of $\slnc$ the subspace
$\mathfrak{a}_k$ consists of matrices with only nonvanishing entries
above the diagonal in the $(N-k+1)$-th column.

We have $\dim\mathfrak{a}_k=N-k$ and
$\np=\mathop{\oplus}\limits_{k=1}^{N-1}\mathfrak{a}_{k}$.
 Each $\mathfrak{a}_k$
is an abelian subalgebra of $\mathfrak{b}_+$ and $\mathfrak{a}_1$ is
an abelian ideal of $\mathfrak{b}_+$. Let us also define
\begin{equation}\label{bk}
\mathfrak{b}_k=\bigoplus_{l=k}^{N-1}\mathfrak{a}_l.
\end{equation}
Obviously
$\mathfrak{a}_k\subset\mathfrak{b}_k\subset\mathfrak{b}_{k-1}$.
Moreover, each $\mathfrak{b}_k$ is a subalgebra of $\np$ and
$\mathfrak{a}_k$ is an abelian ideal of $\bp_k$, i.e.,
\begin{equation}\label{abcomut}
[\mathfrak{b}_k,\mathfrak{b}_k]\subset\mathfrak{b}_k,\quad
[\mathfrak{a}_k,\mathfrak{a}_k]=0, \quad
[\mathfrak{a}_k,\mathfrak{b}_k]\subset\mathfrak{a}_k.
\end{equation}
The subalgebras of $\nm$ will be denoted by $\am_k$. Each $\am_k$ is the
Hermitian conjugate of $\ap_k$ and is generated by the basis elements
(\ref{Bl}) in the following way:
\begin{eqnarray}\label{setideals_m}
\am_k=\spn\{X_{j_k},X_{j_k+1}\ldots,X_{j_k+N-k-1}\}, \quad
j_k=N(N-k)-\frac{k(k+1)}{2}.
\end{eqnarray}
The subalgebras $\bm_k$ are defined in analogy to (\ref{bk}) and together wit
$\am_k$ they follow relations analogous to (\ref{abcomut}). We have the
following decomposition of $\g=\slnc$ into semidirect sum of $2N-1$ commuting
subalgebras
\begin{equation}\label{commdecomp}
\g=\h\oplus\bigoplus_{i=1}^{N-1}\ap_i\oplus\bigoplus_{j=1}^{N-1}\am_j.
\end{equation}

\section{Properties of the adjoint endomorphism}
\label{sec:adjoint}

Expression (\ref{eq-wn0}) consists of products of exponents of
$\ad{X_k}$. In this section we present some fundamental properties
of the adjoint endomorphisms corresponding to root vectors in
$\slnc$ which are crucial for usefulness of the Wei-Norman method in
case of $SL(N,\C)$ group.
\begin{Lemma}\label{nilpot}
Let $\alpha$ be the root and $X_\alpha\in\slnc$ be the corresponding
root vector. Then:
\begin{enumerate}
\item the image of
$\left(\ad{X_\alpha}\right)^2$ is equal to $\g_\alpha$  and \\
$\ker\left(\left(\ad{X_\alpha}\right)^2\right)=
\spn\left\{\{X_1,\ldots,X_{N^2-1}\}\setminus\{X_{-\alpha}\}\right\}.$
\item $\left(\ad{X_\alpha}\right)^3=0$.
\end{enumerate}
\end{Lemma}

\proof First observe that statement 2 follows from 1. To prove the latter
recall that \cite{humphreys}
\begin{equation}
\ad{X_\alpha}(X_\beta)=\left[X_\alpha,X_\beta\right]=
\begin{cases}
H_\alpha, & \alpha+\beta=0,\\
N_{\alpha,\beta}X_{\alpha+\beta}, & \alpha +\beta\in\Delta ,\\
0, & \mbox{otherwise},
\end{cases}
\end{equation}
where $H_\alpha\in \h$ corresponds to root $\alpha$ and
$N_{\alpha,\beta}$ is some constant. We have
\begin{equation}\label{adx2}
\left(\ad{X_\alpha}\right)^2(X_\beta)=\left[X_\alpha,
\left[X_\alpha,X_\beta\right]\right]=
\begin{cases}
-\alpha(H_\alpha)X_\alpha, & \alpha+\beta=0,\\
N_{\alpha,\beta}N_{\alpha,\alpha+\beta}
X_{\alpha+\alpha+\beta}, & \alpha +\alpha +\beta\in\Delta ,\\
0, & \mbox{otherwise}.
\end{cases}
\end{equation}
Since the condition $\alpha +\alpha +\beta\in\Delta$ is never fulfilled for
$\slnc$, equation (\ref{adx2}) implies that $\left(\ad{X_\alpha}\right)^2$
sends any element of $\np\cup\nm$ into $\g_\alpha$. And since for $H\in\h$
\begin{equation}
\left(\ad{X_\alpha}\right)^2(H)=[X_\alpha,[X_\alpha,H]]=-\alpha(H)[X_\alpha,X_\alpha]=0
\end{equation}
the only element of basis (\ref{Bu}-\ref{Bl}) on which
$\left(\ad{X_\alpha}\right)^2$ takes the nonzero value is
$X_{-\alpha}$, what ends the proof of statement 1. \qed

Observe that if we have two commuting matrices of given nilpotency
order $r$ then the sum of the matrices is also nilpotent of order
$r$. It follows from the Jordan theorem \cite{humphreys} -
the matrices can be
expressed in Jordan form in the same basis and are block-diagonal
with the same blocks of maximal size $r-1$. Since $\ap_k$ and
$\am_k$ are commuting subalgebras, the lemma \ref{nilpot} yields:

\begin{Corollary}\label{anilpot}
If $X\in\ap_k$ or $X\in\am_k$ then $(\ad{X})^3=0$.
\end{Corollary}

\begin{Lemma}\label{triang}
Let $\alpha$ be the root and $X_\alpha\in\slnc$ be the corresponding
root vector. In the basis defined by (\ref{Bu}-\ref{Bl}) the
matrices of $\ad{X_\alpha}$ for $\alpha\in\Delta_+$ are strictly
upper triangular and for $\alpha\in\Delta_-$ are strictly lower
triangular.
\end{Lemma}
\proof Let $X_\alpha\in\np$. In the previous section we have mentioned the
theorem of Lie, which has a direct consequence that the matrix of
$\ad{X_\alpha}:\np\to\np$ is strictly upper triangular. So the sector of the
matrix of endomorphism $\ad{X_\alpha}:\g\to\g$ corresponding do $\np$ is
strictly upper triangular. Since $[\np,\h]\subset\np$ the only nonzero matrix
elements of $\ad{X_\alpha}$ in the sector corresponding to $\h$ lie above
diagonal. Finally we consider the action $\ad{X_\alpha}(X_\beta)$ for
$X_\beta\in\nm$. The result is nonzero in two cases: $\alpha+\beta=0$ or
$\alpha+\beta\in\Delta$. In the first case $\ad{X_\alpha}(X_\beta)\in\h$ and
the corresponding matrix element is above diagonal. In te second case it may be
easily checked directly from the definition (\ref{ord_np}-\ref{ord_nm}), that
$\alpha+\beta\prec\beta$. This finishes the proof for $X_\alpha\in\np$. The
proof for $X_\beta\in\nm$ is analogous. \qed

\begin{Lemma}\label{invsp}
Let $X_\alpha\in\ap_k\subset\slnc$ or
$X_\alpha\in\am_k\subset\slnc$, where $\alpha$ is the corresponding
root. The subalgebras $\ap_l$, $\am_l$ for $l<k$ and the subalgebra
$\bp_k\oplus\h\oplus\bm_k$ are invariant subspaces of
$\ad{X_\alpha}$.
\end{Lemma}
\proof Let $X_\alpha\in\ap_k$. We consider three cases:
\begin{enumerate}
\item $X_\beta\in\ap_l$, $l<k$. In this case $X_\alpha\in\bp_l$ and
$X_\beta\in\bp_l$, because $\ap_k\subset\bp_k\subset\bp_l$. Thus,
$\ad{X_\alpha}(X_\beta)=[X_\alpha,X_\beta]\in\bp_l$. On the other
hand $X_\beta\in\ap_l$ and $\ap_l$ is an ideal in $\bp_l$, so
$\ad{X_\alpha}(X_\beta)\in\ap_l$.
\item $Y\in\bp_k\oplus\h\oplus\bm_k$. We have also
$X_\alpha\in\bp_k\oplus\h\oplus\bm_k$ and the property
$\ad{X_\alpha}(Y)\in\bp_k\oplus\h\oplus\bm_k$ follows from the
fact that $\bp_k\oplus\h\oplus\bm_k$ is subalgebra of $\g$ and
this follows directly from the definition of $\bp_k$ and
$\bm_k$ (see (\ref{bk})).
\item $X_\beta\in\am_l$, $l<k$. According to definitions
(\ref{setideals}) and (\ref{setideals_m}) the roots are
    \begin{align}
    \alpha&=(i,N-k+1),\qquad i<N-k+1,\\
    \beta&=(N-l+1,j), \qquad N-l+1>j.
    \end{align}
    We have
    \begin{equation}
    \ad{X_\alpha}(X_\beta)=[X_\alpha,X_\beta]\neq0
    \Longleftrightarrow \alpha+\beta\in\Delta
    \Longleftrightarrow i=j \mbox{ or } k=l.
    \end{equation}
    Condition $k=l$ contradicts the assumption $k>l$.
    So the only possibility is $i=j$, what implies
    $\alpha+\beta=(N-l+1,N-k+l)$. Since $k>l$, we have
    $N-l+1>N-k+1$, so $\ad{X_\alpha}(X_\beta)=X_{\alpha+\beta}\in\am_l$.
\end{enumerate}
The same reasoning holds for $X_\alpha\in\am_k$.\qed

Since for a given $k$ we have the following decomposition of $g$:
\begin{equation}\label{blocks}
\g=\ap_1\oplus\ldots\oplus\ap_{k-1}\oplus
\underbrace{\left(
\bp_k\oplus\h\oplus\bm_k
\right)}_{\mbox{one block}}
\oplus\am_{k-1}\oplus\ldots\oplus\am_1,
\end{equation}
the lemmas \ref{nilpot}, \ref{triang} and \ref{invsp} yield
\begin{Corollary}\label{blockdiagonal}
Let $X_\alpha\in\ap_k\subset\slnc$ ($X_\alpha\in\am_k\subset\slnc$)
be the root vector. In the basis (\ref{Bu}-\ref{Bl}) the matrix of
endomorphism $\ad{X_\alpha}$ is nilpotent of order 3, strictly upper
triangular (lower triangular) and block diagonal with respect to
decomposition (\ref{blocks}).
\end{Corollary}

\section{Exponential function on commuting nilpotent subalgebras of $\slnc$}
\label{sec:nilpotent}

The corollary \ref{blockdiagonal} implies the following
\begin{Corollary}\label{expprop}
For $X\in\ap_k$ or
$X\in\am_k$ the matrix of $\exp(\ad{X})$ is a quadratic polynomial
in $\ad{X}$. Moreover in the basis (\ref{Bu}-\ref{Bl}) the matrix
of $\exp(\ad{X})$ is upper triangular (lower triangular for
$X\in\am_k$) and block diagonal with respect to decomposition
(\ref{blocks}).
\end{Corollary}

It is thus convenient to rewrite the equation (\ref{eq-wn0}) as a sum of terms
corresponding to commuting subalgebras of decomposition (\ref{commdecomp}). We
define:
\begin{equation}\label{defU}
\up_k:=\prod_{i\in \indp} \exp(u_i\ad{X_i}),
 \qquad
 \um_k:=\prod_{i\in\indm} \exp(u_i\ad{X_i}),
\end{equation}
where $\indp=\{i_k,\ldots,i_k+N-k-1\}$ is the index range defined in
(\ref{setideals}), so for given $k$ index $i\in \indp$ numerates
generators of subalgebra $\ap_k$, and
$\indm=\{j_k,\ldots,j_k+N-k-1\}$ is the index range defined in
(\ref{setideals_m}), so for given $k$ index $i\in \indm$ numerates
generators of subalgebra $\am_k$. For $l<N-k-1$ we also define:
\begin{equation}\label{defV}
\vp_{kl}:=\prod_{i=i_k}^{i_k+l} \exp(u_i\ad{X_i}),
 \qquad
 \vm_{kl}:=\prod_{i=j_k}^{j_k+l} \exp(u_i\ad{X_i}),
\end{equation}
where for given $k$ index $i$ runs over first $l-1$ elements of $\indp$ or $\indm$.
We set
$$
\up_0=\um_N=\vp_{k0}=\vm_{k0}=\mathbbm{1}.
$$
Observe,
that commutativity of subalgebras $\ap_k$ implies that:
\begin{equation}
\up_k=\prod_{i\in J_k} \exp(u_i\ad{X_i})
 =\exp\left(\sum_{i\in J_k} u_i\ad{X_i} \right)
 =\exp\left(\ad{\left(\sum_{i\in J_k} u_iX_i\right)} \right),
\end{equation}
so $\up_k$ equals $\exp(\ad{X})$ for some $X\in\ap_k$ (analogously
$\um_k$ equals $\exp(\ad{X})$ for some $X\in\am_k$), so $\up_k$ and
$\um_k$ have the triangularity and block diagonality properties
implied by corollary \ref{expprop}  and the same statement holds for
$\vp_{kl}$ and $\vm_{kl}$. From lemma \ref{nilpot} we have also the following crucial property:
\begin{Corollary}\label{quadratic}
The terms quadratic in parameters $u_i$ appear only when $\up_k$ or $\vp_{kl}$ acts on element of $\am_k$ for the same index $k$ and the terms quadratic in parameters $u_i$ appear only when $\um_k$ or $\vm_{kl}$ acts on element of $\ap_k$ for the same index $k$.
\end{Corollary}

Let $\ih=\frac{1}{2}N(N-1)+1$ be the number of the first element of
Cartan subalgebra according to ordering (\ref{Bu}-\ref{Bl}). We
define also:
\begin{equation}\label{defH}
H:=\prod_{i=\ih}^{\ih+N-1}\exp(u_i\ad{X_i}),
 \qquad
 H_l:=\prod_{i=\ih}^{\ih+l}\exp(u_i\ad{X_i}),
\end{equation}
where $l\in\{1,\ldots,N-2\}$ and $X_i$ generate Cartan subalgebra
$\h$ (see (\ref{Bcartan})). $H$ and $H_l$ are diagonal matrices and
$H_0=\mathbbm{1}$.

\section{Wei--Norman method in a properly chosen basis}
\label{sec:basis}

The equation (\ref{eq-wn0}) can be now rewritten in terms of
operators defined in (\ref{defU}), (\ref{defV}) and (\ref{defH}) in
the following simplified way:
\begin{align}\label{msimpl1}
M(t)=&\sum_{k=1}^{N-1}\prod_{l=0}^{k-1}\up_l\cdot\sum_{i=0}^{N-k-2} \vp_{ki}
 \,u^\prime_{i_k+i}X_{i_k+i}+\\ \label{msimpl2}
&+\prod_{j=1}^{N-1}\up_j\cdot\sum_{i=0}^{N-2}H_i
 \,u^\prime_{\ih+i}X_{\ih+i}+\\ \label{msimpl3}
&+\prod_{j=1}^{N-1}\up_j \cdot H \cdot
 \sum_{k=1}^{N-1}\prod_{l=N}^{N-k+1}\um_l\cdot\sum_{i=0}^{N-k-2} \vm_{ki}
 \,u^\prime_{j_k+i}X_{j_k+i}.
\end{align}

On the other hand the matrix $M(t)$ by definition (\ref{Mdec}) is a linear combination of
generators $X_k$ with given coefficients $a_k(t)$. The
first step in order to solve the system of differential equations is
to solve linear algebraic equations for $u^\prime_k$. This is done
by matrix inversion (\ref{eq-wn2a}-\ref{eq-wn3}). We will show now,
that the expression (\ref{msimpl1}-\ref{msimpl3}) allows for
separation of the full set of equations into sectors corresponding
to decomposition (\ref{commdecomp}).

The first term in (\ref{msimpl1}) corresponding to $k=1$ is an
element of $\ap_1$. It follows from corollary \ref{expprop} and the
fact that it is a sum of terms of the form $\ad{X}Y$ for $X\in\ap_1$
and $Y\in\ap_1$. We move this term to the left hand side of equation (\ref{msimpl1}-\ref{msimpl3}). The rest of the sum, remaining on the right hand side has a common factor $\up_1$, so we multiply both sides by $\up_1^{-1}$. After this operations right hand side does not contain terms proportional to elements of
$\ap_1$, because is an composition of action of block diagonal
operators with respect to decomposition (\ref{blocks}) for $k=1$
acting on generators of $\g$ starting from $\ap_2$. The left hand side reads
\begin{equation}\label{umt}
\up_1^{-1}\left(M(t)-\sum_{i=1}^{N-1} \vp_{ki-1}
 \,u^\prime_{i}X_{i}\right)=
 \up_1^{-1}\left(\sum_{j=1}^na_j(t)X_j-\sum_{i=1}^{N-1} \vp_{ki-1}
 \,u^\prime_{i}X_{i}\right),
\end{equation}
where we used the definition (\ref{Mdec}). Since there are no elements spanning
$\ap_1$ on the right hand side, the first $N-1$ components of (\ref{umt}) have
to vanish. These $N-1$ equations depend only on $u_i$ for $i<N$ (since $\up_1$
depends only on them) and $u_i^\prime$ for $i<N$ (since $\vp_{kl}$ are upper
triangular, what follows from corollary \ref{expprop}). Moreover $\up_1$ is a
quadratic function of function $u_i$, so we end up with matrix Riccati equation
for functions $u_i$ for $i<N$. Once this system of equations is solved and its
solutions are substituted into equation (\ref{msimpl1}-\ref{msimpl3}), the
terms corresponding to generators of $\ap_1$ cancel, and we obtain reduced
system of equations in smaller subspace of $\g$. Next we apply the same
procedure to the remaining sum on the right hand side of
(\ref{msimpl1}-\ref{msimpl3}). It starts now from term spanned by generators of
$\ap_2$. The same reasoning as for $\ap_1$ applies and after moving the term in
question to the left hand side and multiplying both sides by $\up_2$ we obtain
a matrix Riccati equations for functions $u_i$ which multiply generators of
$\ap_2$ in (\ref{coord2}). We substitute the solutions of this system of
equations and the main equation reduces again. We keep repeating the procedure
described above until we obtain the solutions for all coefficient functions
$u_i$ corresponding to all $\ap_k$. The solutions will always come from Riccati
type equations. What remains from the sum (\ref{msimpl1}-\ref{msimpl3}) on the
right hand side of (\ref{eq-wn2a}) after inserting the solutions and
multiplying both sides by product of all $\up_i^{-1}$ is the following:
\begin{align}\label{msimplsimpl}
\mathrm{(l.h.s.)}=&\sum_{i=0}^{N-2}H_i
 \,u^\prime_{\ih+i}X_{\ih+i}+ H \cdot
 \sum_{k=1}^{N-1}\prod_{l=N}^{N-k+1}\um_l\cdot\sum_{i=0}^{N-k-2} \vm_{ki}
 \,u^\prime_{j_k+i}X_{j_k+i}.
\end{align}
First sum corresponds to Cartan subalgebra $\h$. Since $\h$ is commutative we have
\begin{equation}
\exp\left(\ad{X}\right)Y=Y, \quad \mbox{for } X,Y\in\h,
\end{equation}
so in (\ref{msimplsimpl}) $H_i X_{\ih+i}=X_{\ih+i}$, what follows from
definition (\ref{defH}). The second sum in  (\ref{msimplsimpl}) is a
combination of generators of subalgebras $\am_k$ only, because is an action of
diagonal operator $H$ and lower triangular operators $\um_l$ and $\vm_{ki}$ on
those generators. Thus the functions $u_i$ corresponding to generators of $\h$
can be found by simple integration of solutions to previously found Riccati
equations. After substitution of these solutions and multiplication of both
sides of (\ref{msimplsimpl}) by $H^{-1}$ we are left with equation where on
both sides there are only terms proportional to generators $X_j$ spanning
$\am_k$, because all other terms have canceled. In order to find the solutions
for remaining functions $u_i$ we proceed as follows. First we multiply both
sides by $\um_1^{-1}$ which belongs on functions $u_i$ corresponding do $\am_1$
only. Derivatives $u_i^\prime$ on the right hand side appear only in the first
term of the sum, because of the block diagonality of $\vm_{kj}$ operators. Thus
this terms separates form the rest and has to be equal to the action of
$\um_1^{-1}$ on the left hand side. The fact that there are no generators of
$\ap_k$ in the equation and corollary \ref{quadratic} imply that this set of
equations will be linear in functions $u_i$. Once we solve it, the solutions
may be substituted into equation, and terms proportional to generators of
$\am_1$ will cancel. Then we multiply both sides by $\um_2$ and proceed in the
same way obtaining the set of linear equations for functions $u_i$
corresponding to $\am_2$. We keep repeating the described procedure until we
find linear equations for all remaining functions $u_i$.

The above described procedure enables the conversion of highly nonlinear
differential equation (\ref{precession}) into hierarchy of matrix Riccati
equations and linear matrix differential equations. This procedure is effective
in the sense, that it provides an algorithm that may be applied directly. The
authors have written a program in \emph{Maple} which performs this procedure
for any given $N$ and tested its successful performance up to $N=10$.

It is worth mentioning that the crucial ingredient for realizing the
described algorithm in practice is the order of the generators
(\ref{Bu}-\ref{Bl}). Once this base is used for computations and the inverse
in (\ref{eq-wn3}) is successfully computed, the separation of the system of
equations comes up automatically. For large $N$ computation of the inverse is
the part of the algorithm with the largest computational complexity, which
scales with $N$ as $N!$. It follows from the fact, that the matrix $A$ in
(\ref{eq-wn2a}) which is in principal of dimension $(N^2-1) \times (N^2-1)$
has a special block upper triangular form with the largest block to invert of
the size $(N-1)\times(N-1)$. The inversion has to be realized by Cramers rule
which has the mentioned computational complexity.


\section{Examples}
\label{sec:examples}

\subsection{$SL(2,\mathbb{C})$}

For $N=2$ the system (\ref{eq-wn3}) reduces to one Riccati equation,
\begin{equation}\label{2riccati}
u_1^\prime=a_1+2a_2 u_1-a_3 u_1^2,
\end{equation}
and two equations for the remaining unknowns,
\begin{eqnarray}\label{2int1}
u_2^\prime&=&a_2-a_3 u_1, \\
u_3^\prime&=&a_3{\rm e}^{2u_2},
\end{eqnarray}
which give $u_2$ and $u_3$ by simple integrations, once solutions of
(\ref{2riccati}) are known.

\subsection{$SL(3,\mathbb{C})$}
For $N=3$ we obtain from (\ref{eq-wn3}):
\begin{enumerate}
\item A system of two coupled Riccati equations
\begin{eqnarray}\label{3riccati1}
   {u}^\prime_1&=&  a_1+(2a_5-a_4)u_1+a_6u_2
-a_8u_1^2-a_7u_1u_2,
\\
  {u}^\prime_2&=& a_2+a_3u_1+(a_4+a_5)u_2
-a_8u_1u_2-a_7u_2^2,
\end{eqnarray}
which for further reference we will rewrite in the form
\begin{equation}\label{3riccati1m}
\mathbf{u}_{(1)}^\prime=\mathbf{c}_{(1)}^{\phantom T}+
C_{(1)}^{\phantom T}\mathbf{u}_{(1)}^{\phantom T}+
\mathbf{u}_{(1)}^{\phantom T}\mathbf{u}_{(1)}^T\mathbf{b}_{(1)}^{\phantom T},
\end{equation}
with
\begin{equation}\label{3riccati1md1}
\mathbf{u}_{(1)}^{\phantom T}=
\left[
\begin{array}{c}
 u_1 \\
 u_2 \\
\end{array}
\right], \quad
\mathbf{c}_{(1)}^{\phantom T}=
\left[
\begin{array}{c}
 a_1 \\
 a_2 \\
\end{array}
\right],  \quad
\mathbf{b}_{(1)}^{\phantom T}=
\left[
\begin{array}{c}
\begin{array}{c}
 -a_8 \\
 -a_7 \\
\end{array}
\end{array}
\right],
\end{equation}
and
\begin{equation}\label{3riccati1md2}
C_{(1)}^{\phantom T}=
\left[
\begin{array}{cc}
 2a_5-a_4 &   a_6   \\
   a_3    & a_4+a_5 \\
\end{array}
\right].
\end{equation}
\item An equation for $u_3$ which reduces to a scalar Riccati equation
    upon substituting solutions of (\ref{3riccati2})
\begin{equation}\label{3riccati2}
u_3^\prime=(a_3-a_8u_2)+(2a_4-a_5+a_8u_1-a_7u_2)\,u_3+(a_7u_1-a_6)\,u_3^2,
\end{equation}

\item Equations for the coefficients of the Cartan algebra generators
    which are solved by single integrations once solutions of
    (\ref{3riccati1}) and (\ref{3riccati2}) are known,
\begin{eqnarray}\label{3cartan}
u_4^\prime&=&a_4-a_6u_3+ a_7\left( u_1u_3-u_2 \right) \\
u_5^\prime&=&a_5-a_8u_1-a_7u_2
\end{eqnarray}

\item Equations for the coefficients of the generators of the second
    ("lower-triangular") nilpotent subalgebra
\begin{eqnarray}
u_6^\prime &=& \left(a_6-a_7u_1\right){\rm e}^{2\,u_4-u_5}, \label{3Bl1} \\
u_7^\prime &=&(a_7u_3+a_8)u_6{\rm e}^{-u_4+2\,u_5} +a_7{\rm
e}^{u_4+u_5}, \label{3Bl2} \\
u_8^\prime &=& (a_8+a_7u_3){\rm e}^{-u_4+2\,u_5}, \label{3Bl3}
\end{eqnarray}
which are solved by simple consecutive integrations.

\end{enumerate}

\subsection{$SL(4,\mathbb{C})$}
A similar structure emerges for $N=4$. The system of $15$ equations
(\ref{eq-wn3}) separates into:
\begin{enumerate}
\item A system three coupled Riccati equations,
\begin{equation}\label{4riccati1}
\mathbf{u}_{(1)}^\prime=\mathbf{c}_{(1)}^{\phantom T}+
C_{(1)}^{\phantom T}\mathbf{u}_{(1)}^{\phantom T}+
\mathbf{u}_{(1)}^{\phantom T}\mathbf{u}_{(1)}^T\mathbf{b}_{(1)}^{\phantom T},
\end{equation}
where
\begin{equation}\label{4riccati1d}
\mathbf{u}_{(1)}^{\phantom T}=
\left[
\begin{array}{c}
 u_1 \\
 u_2 \\
 u_3 \\
\end{array}
\right], \quad
\mathbf{c}_{(1)}^{\phantom T}=
\left[
\begin{array}{c}
 a_1 \\
 a_2 \\
 a_3 \\
\end{array}
\right], \quad
\mathbf{b}_{(1)}^{\phantom T}=
\left[
\begin{array}{c}
 -a_{15} \\
 -a_{14} \\
 -a_{13} \\
\end{array}
\right],
\end{equation}
and
\begin{equation}\label{4riccati1d1}
C_{(1)}^{\phantom T}=
\left[
\begin{array}{ccc}
 -a_8+2a_9 &    a_{12}    & a_{11}  \\
   a_4    & -a_7+a_8+a_9 & a_{10}  \\
   a_5    &     a_6      & a_7+a_9 \\
\end{array}
\right].
\end{equation}

\item A system of two coupled Riccati equations,
\begin{equation}\label{4riccati2}
\mathbf{u}_{(2)}^\prime=\mathbf{c}_{(2)}^{\phantom T}+
C_{(2)}^{\phantom T}\mathbf{u}_{(2)}^{\phantom T}+
\mathbf{u}_{(2)}^{\phantom T}\mathbf{u}_{(2)}^T\mathbf{b}_{(2)}^{\phantom T},
\end{equation}
where
\begin{equation}\label{4riccati2d1}
\mathbf{u}_{(2)}^{\phantom T}=
\left[
\begin{array}{c}
 u_4 \\
 u_5 \\
\end{array}
\right], \quad
\mathbf{c}_{(2)}^{\phantom T}=
\left[
\begin{array}{c}
 a_4-a_{15}u_2 \\
 a_5-a_{15}u_3 \\
\end{array}
\right],  \quad
\mathbf{b}_{(2)}^{\phantom T}=
\left[
\begin{array}{c}
-a_{12}+a_{14}u_1 \\
-a_{11}+a_{13}u_1 \\
\end{array}
\right],
\end{equation}
and
\begin{equation}\label{4riccati2d2}
C_{(2)}^{\phantom T}=
\left[
\begin{array}{cc}
 -a_7+2a_8-a_9-a_{14}u_2+a_{15}u_1 &        a_{10}-u_2a_{13}-a_{15}        \\
           a_6-a_{14}u_3           & a_7+a_8-a_9-a_{13}u_{3}+a_{15}u_1 \\
\end{array}
\right].
\end{equation}
Observe that once a solution of (\ref{4riccati1}) is known, the system
(\ref{4riccati2}) is closed since (\ref{4riccati2d1}) and
(\ref{4riccati2d2}) are given in terms of some known functions.

\item A scalar Riccati equation,
\begin{eqnarray}\label{4riccati3}
u_6^\prime&=&a_{{6}}-a_{{12}}u_{{5}}+a_{{14}}u_{{1}}u_{{5}}-a_{{14}}u_{{3}}
\nonumber \\
&+&(2\,a_{{7}}-a_{{8}}-u_{{5}}a_{{11}}+u_{{4}}a_{{12}}-u_{{3}}a_{{13}}+a_{
{13}}u_{{1}}u_{{5}}-a_{{14}}u_{{1}}u_{{4}}+u_{{2}}a_{{14}}
)u_6
\nonumber \\
&+&(-a_{{10}}+u_{{4}}a_{{11}}-a_{{13}}u_{{1}}u_{{4}}+a_{{13}}u_{{2}})u_6^2,
\end{eqnarray}
with coefficients depending on solutions of (\ref{4riccati1}) and
(\ref{4riccati2}).
\item The remaining 9 equations for $u_7,\ldots,u_{15}$ which are solved
    by single interactions of functions constructed from the initial
    coefficients $a_7,\ldots,a_{15}$, and solutions of (\ref{4riccati1}),
    (\ref{4riccati2}) and (\ref{4riccati3}).
\end{enumerate}

\section*{Acknowledgments}

The presented results are obtained in frames of the the Polish National Science
Center project MAESTRO DEC-2011/02/A/ST1/00208 support of which is gratefully
acknowledged by both authors.


\begin{thebibliography}{10}

\bibitem{wei63} J.~Wei and E.~Norman.
\newblock Lie algebraic solution of linear differential equations.
\newblock {\em J. Math. Phys.}, 4(4):575--581, 1963.

\bibitem{wei64} J.~Wei and E.~Norman.
\newblock On global representations of the solutions of linear differential
  equations as a product of exponentials.
\newblock {\em Proc. Am. Math. Soc.}, 15(2):327--334, 1964.

\bibitem{brockett72} R.~Brockett.
\newblock Systems theory on groupmanifolds and coset spaces.
\newblock {\em SIAM J. Control}, 10:265--284, 1972.

\bibitem{brockett73} R.~Brockett.
\newblock Lie algebras and lie groups in control theory.
\newblock In R.~Brockett D.~Maine, editor, {\em Geometric Methods in Systems
  Theory, Proceedings of the NATO Advanced Study Institute}, D. Reidel, Dordrecht, 1973.

\bibitem{huillet87} T.~Huillet, A.~Monin, and G.~Salut.
\newblock Minimal realizations of the matrix transition lie groupfor bilinear
  control systems: explicit results.
\newblock {\em Systems Control Lett.}, 9:267–--274, 1987.

\bibitem{leonard93} N.~Leonard and P.~Krishnaprasad.
\newblock Averaging on lie groups, attitude control and drift.
\newblock In {\em Proceedings of the Twenty-Seventh Annual Conference on
  Information Sciences and Systems}, 369--374, 1993.

\bibitem{chiou94} W.~Chiou and S.~Yau.
\newblock Finite dimensional filters with nonlinear drift, brocketts problem on
  classification of finite dimensional estimation algebras.
\newblock {\em SIAM J. Control Optim.}, 32:297–--310, 1994.

\bibitem{leonard95} N.~Leonard and P.~Krishnaprasad.
\newblock Motion control of drift-free left-invariant systems on lie groups.
\newblock {\em IEEE Trans. Automat. Control}, 40:1539–--1554, 1995.

\bibitem{charalambous00} C.D. Charalambous and R.J. Elliott.
\newblock Information states in stochastic control and filtering: A lie
  algebraic theoretic approach.
\newblock {\em IEEE Trans. Automat. Control}, 45(4):653--674, 2000.

\bibitem{carinena03} J. Cari{\~n}ena, J. Clemente-Gallardo, and
    Arturo Ramos.
\newblock Motion on {L}ie groups and its applications in control theory.
\newblock {\em Rep. Math. Phys.}, 51:159--170, 2003.

\bibitem{owren01} B.~Owren and A.~Marthinsen.
\newblock Integration methods based on canonical coordinates of the second
  kind.
\newblock {\em Numer. Math.}, 87:763–-790, 2001.

\bibitem{celledoni01} E.~Celledoni and A.~Iserles.
\newblock Methods for the approximation of the matrix exponential in a
  {L}ie-algebraic setting.
\newblock {\em J. Numer Anal.}, 21:463--488, 2001.

\bibitem{zanna02} A.~Zanna and H.~Z. Munthe-Kaas.
\newblock Generalized polar decompositions for the approximation of the matrix
  exponential.
\newblock {\em SIAM J. Matrix Anal.}, 23:840–--862, 2002.

\bibitem{dattoli87} G.~Dattoli, P.~Di~Lazzaro, and A.~Torre.
\newblock $SU(1,1)$, $SU(2)$, and $SU(3)$ coherence-preserving {H}amiltonians and
  time-ordering techniques.
\newblock {\em Phys. Rev. A}, 35:1582--1589, 1987.

\bibitem{carinena02} J.~F.~Cari{\~n}ena and A. Ramos.
\newblock A new geometric approach to lie systems and physical applications.
\newblock {\em Acta Applicandae Mathematicae}, 70:43--69, 2002.

\bibitem{carinena03a} J.~F. Cari{\~n}ena and A.~Ramos.
\newblock Applications of {L}ie systems in quantum mechanics and control
  theory.
\newblock [in] {\em Classical and quantum integrability: dedicated to
  W{\l}odzimierz Tulczyjew}, vol.~59, p. 143. Institute of Mathematics,
  Polish Academy of Sciences, 2003.

\bibitem{carinena09} J.~F. Cari\~nena, J. de~Lucas, and A. Ramos.
\newblock A geometric approach to time evolution operators of lie quantum
  systems.
\newblock {\em Int. J. Theor. Phys.}, 48:1379--1404,
  2009.

\bibitem{kuna10} M.~Kuna and J.~Naudts.
\newblock General solutions of quantum mechanical equations of motion with
  time-dependent {H}amiltonians: A {L}ie algebraic approach.
\newblock {\em Rep. Math. Phys.}, 65(1):77 -- 108, 2010.

\bibitem{nihtila11} M.~Nihtil\"a.
\newblock Wei-{N}orman {T}echnique for {C}ontrol {D}esign of {B}ilinear
{ODE} {S}ystems with {A}pplication to {Q}uantum {C}ontrol,
\newblock [in] Jean L\'evine and Philippe M\"ullhaupt, editors, {\em Advances in
  the Theory of Control, Signals and Systems with Physical Modeling}, volume
  407 of {\em Lecture Notes in Control and Information Sciences}, pp
  189--199. Springer Berlin/Heidelberg, 2011.

\bibitem{humphreys} J.~E.~Humphreys.
\newblock Introduction to Lie Algebras and Representation Theory.
\newblock Springer-Verlag 1980.

\end{thebibliography}

\end{document}